\documentclass[a4,10pt]{article}

\usepackage{wrapfig}

\usepackage{lineno}

\usepackage{amsmath}
\usepackage{amsthm}
\usepackage{amssymb}
\usepackage{enumerate}
\usepackage[pdftex]{graphicx}
\usepackage{float}

\usepackage{url}
\usepackage[colorlinks=true, linkcolor=blue, citecolor=blue, urlcolor=blue]{hyperref}
\usepackage{cite}

\usepackage[top=36truemm,bottom=35truemm,left=32truemm,right=32truemm]{geometry}

\usepackage{titlesec}
\titleformat*{\section}{\large\bfseries}

\theoremstyle{definition}
\newtheorem{theorem}{Theorem}[section]

\newtheorem{lemma}[theorem]{Lemma}

\newtheorem{proposition}[theorem]{Proposition}
\newtheorem{remark}[theorem]{Remark}
\numberwithin{equation}{section}

\allowdisplaybreaks

\begin{document}
\title{A rigorous Abel-summed Hardy-type identity \\for $p$-circle domains via regularized Poisson summation}
  \author{\large{Masaya Kitajima}}
  \date{}
  \maketitle
  \begin{abstract}
  In our previous work, \textit{a generalized Hardy-type identity} (GHI) for \textit{astroid-type $p$-circle} domains was formally derived in terms of generalized $p$-Bessel functions. However, the direct application of the Poisson summation formula to the indicator functions of these domains leaves the convergence and rigorous justification of the resulting oscillatory series unresolved. In the present paper, we address this difficulty by introducing an anisotropic regularization based on convolution kernels with Fourier multipliers $e^{-t(|\xi_{1}|^p+|\xi_{2}|^{p})}$. We establish the corresponding regularized Poisson summation identity and rigorously justify the limiting process through uniform estimates for lattice-point sums. A key feature of our analysis is the determination of boundary contributions: smooth boundary points carry weight $1/2$, whereas the contributions from axial cusps vanish for $0<p<1$. These results lead to an Abel-summed GHI, expressed in terms of the distorted-angle-dependent $p$-Bessel functions introduced in our previous work, with explicit boundary corrections. In particular, for $p=2$, the identity recovers the \textit{Abel-summed form of the classical Hardy--Vorono\"i Bessel series}. Our result provides a rigorous regularized counterpart of the previously proposed identity, while the convergence of the corresponding unregularized series remains a separate problem.\\
\textbf{Keywords:} Fourier analysis, generalized Hardy-type identity, $p$-Bessel functions, Abel summation, Poisson summation; boundary corrections, $p$-circle domains.\\
\textbf{2020 Mathematics Subject Classification:} Primary 42B05, Secondary 42B08, 33C10, 11P21.
  \end{abstract}\vspace{-15pt}
\begingroup
\renewcommand{\thefootnote}{}
\footnotetext{This work was supported by JSPS KAKENHI Grant Number JP26KJ1294.}
\endgroup
  
  \section{Introduction and main results}
  \hspace{13pt}The classical \textit{Hardy–Vorono\"{i} identity} expresses the lattice-point discrepancy $P(r):=N(r)-\pi r^{2}$ for the Euclidean disk in terms of an oscillatory Bessel series:
\begin{equation*}\label{HI}
  P(r)=r\sum_{k=1}^{\infty}\frac{R(k)}{\sqrt{k}}J_{1}(2\pi r\sqrt{k})\quad\text{for }r^{2}\notin\mathbb{N},
\end{equation*}
where $N(r)$ denotes the number of lattice points in the open disk of radius $r$, and $R(k):=\#\{n\in\mathbb{Z}^{2}\mid|n|^{2}=k\}$. This formula was first proposed by G. Vorono\"{i} in 1905 through arguments that were not fully rigorous, and a rigorous proof was subsequently given by G.H. Hardy \cite{Hardy-1915}. Further developments, including E. Landau's approach and the convergence properties of the associated Bessel series, are discussed in \cite[Chapter I]{Popov}. In particular, Hardy's later work \cite{Hardy-1917} illustrates the careful analysis of convergence and limiting procedures required in deriving related Bessel-series formulas. This historical development highlights the distinction between formally derived oscillatory identities and rigorously justified summation formulas, which also motivates the present study.\par
In our previous work \cite[(2.8)--(2.9)]{K3}, we introduced generalized $p$-Bessel functions depending on distorted angles, defined by the following series representations. \vspace{-3pt}
\begin{equation*}\label{p-Bessel-S}\vspace{-2pt}
\mathcal{J}_{\omega,\varphi}^{[p]}(r)
= \Bigl(\frac{2}{p}\Bigr)^{2+\omega} \frac{\pi}{\Gamma(\frac{1}{p})^2} 
\sum_{k=0}^{\infty} \frac{(-1)^k}{k! \, \Gamma(\frac{2}{p}(k+1) + \omega)}
\Bigl(\frac{r}{2}\Bigr)^{2k+\omega} \Phi_{k,\varphi}^{[p]}, 
\quad r \geq 0,
\end{equation*}
\begin{equation*}
\text{with }\Phi_{k,\varphi}^{[p]} 
= \sum_{n=0}^{k} \frac{\bigl(\frac{2}{p}(k+1)-1\bigr)!}{n!\,(k-n)!} 
\left( \frac{B(\frac{2}{p}(n+\frac{1}{2}), \frac{2}{p}(k-n+\frac{1}{2}))}{B(n+\frac{1}{2}, k-n+\frac{1}{2})} \right) 
(\cos^{\frac{4}{p}}\varphi)^n (\sin^{\frac{4}{p}}\varphi)^{k-n}.
\end{equation*}
We subsequently proposed a \textit{generalized Hardy-type identity} (GHI) for the domains bounded by \textit{astroid-type $p$-circle}, defined by $\{x\in\mathbb{R}^{2}\mid|x_{1}|^{p}+|x_{2}|^{p}\leq r^{p}\}$ with $p > 0$ satisfying $2/p \in \mathbb{N}$ \cite[Theorem 1.2]{K3}. However, the identity was obtained through a formal application of the \textit{Poisson summation formula} to indicator functions, and the convergence of the resulting unregularized oscillatory series remains unresolved. The purpose of the present paper is to establish a rigorously justified \textit{Abel-summed} counterpart of this identity. \par
To this end, we employ anisotropic convolution kernels whose Fourier transforms are given by $\widehat{K_{t}^{[p]}}(\xi_1,\xi_2)=e^{-t(|\xi_1|^p+|\xi_2|^p)}$. Although regularization of the Poisson summation formula is a classical technique, its application to the present setting requires a careful justification of the limiting process. In particular, we determine the contributions of boundary lattice points and show that smooth boundary points carry weight $1/2$, whereas the contributions from axial cusps vanish for $0<p<1$.\par
Combining these boundary contributions with the Fourier representation involving the generalized $p$-Bessel functions from \cite[(2.3)]{K3}, we establish an Abel-summed GHI with explicit boundary corrections. The resulting identity also recovers the Abel-summed form of the classical Hardy–Voronoï Bessel series when $p=2$. Let $\mathcal{A}_{s}^{[p]}$ be the finite set of distorted angles $\varphi$ corresponding to lattice points on the $p$-circle of radius $s^{1/p}\ (\geq1)$, namely, 
\begin{equation*}\label{mathcalA}
  \mathcal{A}_{s}^{[p]}:=\left\{\varphi\in[0,2\pi)\ \middle|\ \left(\mathrm{sgn}\left(\cos\varphi\right)s^{\frac{1}{p}}|\cos\varphi|^{\frac{2}{p}},\mathrm{sgn}\left(\sin\varphi\right)s^{\frac{1}{p}}|\sin\varphi|^{\frac{2}{p}}\right)\in\mathbb{Z}^{2}\right\},
\end{equation*}
and let $\mu$ be the counting measure. Then our main result is stated as follows.

\begin{theorem}[\textit{Abel-summed generalized Hardy-type identity}]\label{thm:abel-GHI}\itshape{
Let $p>0$ satisfy $2/p\in\mathbb{N}$, and let $r>0$. Define $D_{p}(r):=\{x\in\mathbb{R}^{2} \mid |x_{1}|^{p}+|x_{2}|^{p}<r^{p}\}$, and $N_{p}(r):=\#(D_{p}(r)\cap\mathbb{Z}^{2})$. The lattice-point discrepancy $P_{p}$ and the boundary correction $B_{p}$ are defined, respectively, by }
\begin{equation*}
P_{p}(r):=N_{p}(r)-\frac{2\Gamma(1/p)^{2}}{p\Gamma(2/p)}r^{2},\quad B_{p}(r):=\sum\nolimits_{n\in\partial D_{p}(r)\cap\mathbb{Z}^{2}}w_{p}(n)
\end{equation*}
\begin{equation*}
\text{where }w_{p}(n):=
\begin{cases}
1/2, & n_{1}n_{2}\neq 0,\\[2mm]
0, & n_{1}n_{2}=0,\quad 0<p<1,
\end{cases}
\quad
\begin{cases}
1/4, & n_{1}n_{2}=0,\quad p=1,\\[2mm]
1/2, & n_1n_2=0,\quad p=2.
\end{cases}
\end{equation*}
Then
\begin{equation}\label{eq:abel-GHI}
P_{p}(r)=\frac{p\Gamma(1/p)^2}{2\pi}\ r\lim_{t\to0+}\int_{1}^{\infty}\frac{e^{-ts}}{s^{1/p}}\Bigl(\sum\nolimits_{\varphi\in\mathcal{A}_{s}^{[p]}}\mathcal{J}_{1,\varphi}^{[p]}(2\pi s^{\frac{1}{p}}r)\Bigr)\ d\mu(s)-B_{p}(r).
\end{equation}
Moreover, for every $t>0$, this regularized series is absolutely convergent.
\end{theorem}
\begin{remark}[Boundary corrections and the classical case]
An essential feature of the preceding identity is the explicit boundary correction arising from the anisotropic regularization. At smooth boundary points, the limiting contribution is $1/2$, whereas at the axial cusps of the $p$-circle for $0<p<1$, it vanishes. In the case $p=1$, the axial vertices contribute $1/4$. These boundary contributions arise naturally from the limiting behavior of the regularization and are accounted for by the correction term $B_{p}(r)$.\par
In particular, for $p=2$, the distorted-angle dependence disappears, since $\mathcal{J}_{1,\varphi}^{[2]}=J_{1}$. Here, we interpret $R(r^{2})$ as zero when $r^{2}\notin\mathbb{N}$. Then, we obtain 
\begin{equation*}
P_{2}(r)=r\lim_{t\to0+}\sum_{k=1}^{\infty}e^{-tk}\frac{R(k)}{\sqrt{k}}J_{1}(2\pi r\sqrt{k})-\frac{1}{2} R(r^{2}).
\end{equation*}
This formula recovers precisely, in the sense of Abel summation and with the boundary correction corresponding to interior lattice-point counting, the classical Hardy–Vorono\"{i} identity stated in \cite[Section 4, Theorem 4]{Popov}. This reduction provides an important consistency check for the generalized identity, although it does not establish ordinary convergence of the unregularized series for general $p$.
\end{remark}
The remainder of this paper is organized as follows. In Section 2, we construct Gaussian-based smoothing kernels and establish their fundamental properties, including the decay estimates needed for the subsequent analysis. In Section 3, we investigate the uniform behavior of the regularized lattice-point sums and determine their limiting contributions at interior, exterior, and boundary lattice points. In Section 4, we establish the regularized Poisson summation identity and complete the proof of Theorem \ref{thm:abel-GHI}. Finally, Section 5 discusses several open questions and possible directions for future research.

  \section{Construction and Properties of Gaussian-Based Smoothing Kernels for $0<p\leq 2$}
  \hspace{13pt}In this section, we introduce a family of smoothing kernels based on Gaussian functions. We emphasize that these kernels play an important role not only in the subject of the present paper but also in the Fourier analysis of $p$-radial functions. Here, we say that a function $f$ on $\mathbb{R}^{2}$ is \textit{$p$-radial} if there exists a function $\phi_{f}$ on $[0,\infty)$ such that $f(x)=\phi_{f}(|x|_{p})$ for all $x\in\mathbb{R}^{2}$, where $|x|_{p}:=(|x_{1}|^{p}+|x_{2}|^{p})^{\frac{1}{p}}$.\par
  Our construction begins with a suitable family of normalized Gaussian functions and proceeds through a weighted superposition with respect to a probability measure. We then construct a two-dimensional kernel by taking the product of two copies of the resulting one-dimensional kernel. We establish several fundamental properties of these kernels, including positivity, normalization, and concentration at the origin. These properties ensure that the resulting family of kernels forms an approximate identity, providing a useful framework for smoothing functions while preserving the $p$-dependent structure relevant to our analysis. \vspace{12pt}\par

For $u>0$, we define the Gaussian function used throughout this paper by
\begin{equation*}
g_{u}(x):=\sqrt{\frac{\pi}{u}}\ e^{-\frac{\pi^{2}x^{2}}{u}}.
\end{equation*}
Clearly, $g_{u}\geq0$ and $\int_{\mathbb{R}} g_{u}=1$. Moreover, its Fourier transform is given by$\widehat{g_{u}}(\xi)=e^{-u\xi^{2}}$.\par
Indeed, let $a=\pi^{2}/u$ and set $I(\xi):=\int_{\mathbb{R}} e^{-ax^{2}}e^{-2\pi i\xi x}dx$.
Since $|x|e^{-ax^{2}}$ is integrable, differentiation under the integral sign is justified, and hence
\begin{equation*}
\frac{d}{d\xi}I(\xi)=-2\pi i\int_{\mathbb{R}} xe^{-ax^2}e^{-2\pi i\xi x}dx=\frac{\pi i}{a}\int_{\mathbb{R}} \left(\frac{d}{dx}\ e^{-ax^2}\right)e^{-2\pi i\xi x}dx=-\frac{2\pi^2\xi}{a}I(\xi).
\end{equation*}
Thus, $(d/d\xi)e^{\pi^2\xi^2/a}I(\xi)=0$. Since $I(0)=\sqrt{\pi/a}$, we obtain $I(\xi)=\sqrt{\pi/a}e^{-\pi^{2}\xi^{2}/a}$.\vspace{8pt}\par
We now recall the following theorem concerning a probability measure required for our regularization procedure.
\begin{theorem}[\textit{\cite[Chapter 5, Theorem 5.2]{Schilling}}]\label{thm:sub}\itshape{
Let $0<\alpha<1$ and $t>0$. Then there exists a unique probability measure $\nu_{\alpha,t}$ on $[0,\infty)$ such that, for every $\lambda\geq0$,}
\begin{equation*}\label{eq:laplace}
\int_{(0,\infty)}e^{-\lambda u}d\nu_{\alpha,t}(u)=e^{-t\lambda^{\alpha}}.
\end{equation*}
\end{theorem}
For $0<p<2$, let $\alpha=p/2$, and define, for $x\in\mathbb{R}$ and $\eta\in\mathbb{R}^{2}$,
\begin{equation}\label{eq:kernel}
k^{[p]}_{t}(x):=\int_{(0,\infty)}g_{u}(x)d\nu_{\alpha,t}(u),\qquad
K^{[p]}_{t}(\eta):=k^{[p]}_{t}(\eta_{1})\ k^{[p]}_{t}(\eta_{2}).
\end{equation}
Since $k^{[p]}_{t}\geq0$, Tonelli's theorem yields $\int_{\mathbb{R}}k^{[p]}_{t}=1$.
Furthermore, by (\ref{eq:kernel}) and Fubini's theorem,
\begin{equation*}
\widehat{k^{[p]}_{t}}(\xi)=\int_{(0,\infty)}\widehat{g_{u}}(\xi)d\nu_{\alpha,t}(u)=\int_{(0,\infty)} e^{-u\xi^{2}}d\nu_{\alpha,t}(u)=e^{-t|\xi|^{p}}.
\end{equation*}
For $p=2$, on the other hand, we set $k^{[2]}_{t}=g_{t}$. In either case, we have
$\widehat{K^{[p]}_{t}}(\eta)=e^{-t|\eta|_{p}^{p}}$. Moreover,
\begin{equation*}
\int_{\mathbb{R}} e^{-t|\xi|^p} d\xi=2\int_{0}^{\infty}e^{-v}\Bigl(\frac{t^{-\frac{1}{p}}}{p}\Bigr)v^{\frac{1}{p}-1}dv=
\frac{2\Gamma(\frac{1}{p})}{pt^{\frac{1}{p}}}<+\infty.
\end{equation*}
Consequently, by the Fourier inversion theorem, we may choose the continuous representative of the $L^{1}$ equivalence class given by
\begin{equation}\label{eq:inverse}
k^{[p]}_{t}(x)=\int_{\mathbb{R}} e^{-t|\xi|^p}e^{2\pi ix\xi}\ d\xi,
\end{equation}
which satisfies $|k^{[p]}_{t}(x)|\leq(2/p)\Gamma(1/p)t^{-1/p}$. Its continuity follows from the dominated convergence theorem, since the integrand converges pointwise as $x_{j}\to x$ and $e^{-t|\xi|^{p}}\in L^{1}(\mathbb{R})$. The evenness $k^{[p]}_{t}(-x)=k^{[p]}_{t}(x)$ follows from the change of variables $\xi\mapsto-\xi$. Furthermore, if the continuous representative (\ref{eq:inverse}) were negative at some point, it would also be negative in a neighborhood of that point, contradicting the nonnegativity of (\ref{eq:kernel}) and their almost-everywhere agreement. Therefore, the continuous representative is nonnegative at every point.\vspace{8pt}\par
Applying the scaling transformation $\xi\mapsto t^{1/p}\xi$ to (\ref{eq:inverse}), we obtain
\begin{equation}\label{eq:scale}
k^{[p]}_{t}(x)=t^{-1/p}k^{[p]}_{1}(t^{-1/p}x),\qquad
K^{[p]}_{t}(\eta)=t^{-2/p}K^{[p]}_{1}(t^{-1/p}\eta).
\end{equation}
Hence, for any $\varepsilon>0$, the substitution $z=t^{1/p}y\in\mathbb{R}^{2}$ gives
\begin{align}\label{eq:approx}
\int_{|z|_2\ge\varepsilon}K^{[p]}_{t}(z)dz
&=\int_{|y|_2\ge\varepsilon t^{-1/p}}K^{[p]}_{1}(y)dy\longrightarrow0\qquad\text{as }t\to0+.
\end{align}
The last limit follows from $K^{[p]}_{1}\geq0$, $\int_{\mathbb{R}^{2}} K^{[p]}_{1}=1<+\infty$, and the fact that the radius of the integration region satisfies $\varepsilon t^{-1/p}\to\infty$. \\
\begin{proposition}\label{lem:tail}
\itshape{Let $0<p\leq 2$ and $0<\delta<p$. Then there exists a constant $C_{\delta}^{[p]}>0$ such that, for all $t>0$ and $x\neq 0$, the following estimate holds:}
\begin{equation*}\label{eq:tail}
0\leq k^{[p]}_{t}(x)\leq C_{\delta}^{[p]}t^{\delta/p}|x|^{-1-\delta}.
\end{equation*}
\end{proposition}
\begin{proof}
First, suppose that $0<p<2$, and set $q:=\delta/2<\alpha:=p/2<1$. Since $z^{q+1/2}e^{-z}\to0$ as $z\to0+$ and as $z\to\infty$, continuity implies that there exists a constant $C_{q}>0$ such that $e^{-z}\leq C_{q}z^{-q-1/2}$ for all $z>0$. Taking $z=\pi^{2}x^{2}/u$, we obtain
\begin{equation*}
g_{u}(x)=\sqrt{\pi/u}\ e^{-\pi^{2}x^{2}/u}\leq C_{q}' u^{q}|x|^{-1-2q}.
\end{equation*}
Integrating this inequality with respect to $\nu_{\alpha,t}$ and applying the formula \cite[Appendix A, Lemma A.1]{Knopova}
\begin{equation*}\label{eq:moment}
\int_0^\infty u^q\,d\nu_{\alpha,t}(u)
=t^{q/\alpha}\frac{\Gamma(1-q/\alpha)}{\Gamma(1-q)} \quad\text{with }0<q<\alpha<1,
\end{equation*}
we obtain the desired estimate, since $2q=\delta$ and $q/\alpha=\delta/p$.
More precisely,
\begin{equation*}
k^{[p]}_{t}(x)\leq C_{q}'|x|^{-1-2q}t^{q/\alpha}
\frac{\Gamma(1-q/\alpha)}{\Gamma(1-q)}\leq \left(\frac{C_{\delta/2}'\Gamma(1-\delta/p)}{\Gamma(1-\delta/2)}\right)t^{\delta/p}|x|^{-\delta-1}.
\end{equation*}
For $p=2$, the same estimate follows by applying the above bound for $e^{-z}$ directly to $k^{[2]}_{t}=g_{t}$.
\end{proof}

\section{Uniform Estimates at Distant Lattice Points and Regularization Weights}
\hspace{13pt}Let $F_{r}:=\mathbf{1}_{D_{p}(r)}$ be the indicator function of the $p$-circle domain $D_{p}(r)$, and let $K^{[p]}_{t}$ be the anisotropic smoothing kernel constructed in Section 2. For $t>0$ and $n\in\mathbb{Z}^{2}$, consider the regularized lattice-point contribution
\begin{equation*}
G_{t}(n):=(F_{r}\ast K^{[p]}_{t})(n)=\int_{\mathbb{R}^{2}}F_{r}(n-y)K^{[p]}_{t}(y)dy.
\end{equation*}
Since $K^{[p]}_{t}$ is nonnegative and has total integral one, $G_{t}(n)$ represents the proportion of the kernel mass centered at $n$ that falls inside $D_{p}(r)$. In particular, $0\leq G_{t}(n)\leq1$. \par
The purpose of this section is to justify the limiting behavior of the regularized lattice-point sum as $t\to 0+$. To this end, two issues must be addressed. First, we establish a uniform estimate for the contributions from lattice points sufficiently far from $D_{p}(r)$, which allows us to control the tail of the infinite sum. We then determine the pointwise limits of $G_{t}(n)$ at interior, exterior, and boundary lattice points. The boundary limits depend on the local geometry of the $p$-circle and give rise to the correction term $B_{p}(r)$ in our main theorem.\vspace{8pt}\par

We begin with the following estimate for distant lattice points.
\begin{lemma}[\textit{Uniform estimates at distant lattice points}]\label{lem:lattice}
\itshape{Let $0<p\leq 2$, $0<\delta<p$, $R>2r+2$, and $t>0$. 
Then there exists a constant $C_{r,\delta}^{[p]}$ such that}
\begin{equation}\label{eq:lattice}
\sum_{\substack{n\in\mathbb{Z}^2\\|n|_{\infty}>R}}(F_{r}\ast K^{[p]}_{t})(n)\leq C_{r,\delta}^{[p]}t^{\delta/p}R^{-\delta}.
\end{equation}
In particular, the lattice point series converges uniformly with respect to $t\in(0,1]$.
\begin{equation}\label{lim:sup}
\sup_{0<t\leq1}\sum_{|n|_{\infty}>R}(F_{r}\ast K^{[p]}_{t})(n)\leq C_{r,\delta}^{[p]}R^{-\delta}\longrightarrow0\qquad\text{as }R\to\infty.
\end{equation}
\end{lemma}
\begin{proof}
Since $D_{p}(r)\subset[-r,r]^{2}$ and $K^{[p]}_{t}\geq0$, we have
\begin{equation*}
0\leq(F_{r}\ast K^{[p]}_{t})(n)=\int_{D_{p}(r)}K^{[p]}_{t}(n-y)dy\leq\prod_{j=1}^{2}\left(\int_{-r}^{r} k^{[p]}_{t}(n_{j}-y_{j})dy_{j}\right).
\end{equation*}
If $|n|_{\infty}>R$, then either $|n_{1}|>R$ or $|n_{2}|>R$. We first consider the former case. By Tonelli's theorem, we may rearrange the nonnegative sums and separate the resulting expression into the product of two series, which will be estimated individually:
\begin{equation*}
\sum_{|n_{1}|>R}\sum_{n_{2}\in\mathbb{Z}}(F_{r}\ast K^{[p]}_{t})(n)\leq\biggl(\sum_{|n_{1}|>R}\int_{-r}^{r}k^{[p]}_{t}(n_{1}-y_{1})dy_{1}\biggr)\times\biggl(\sum_{n_{2}\in\mathbb{Z}}\int_{-r}^{r}k^{[p]}_{t}(n_{2}-y_{2})dy_{2}\biggr).
\end{equation*}
Since $|n_{1}|>R>2r+2$, we have $-r>-R/2+1>-|n_{1}|/2+1$. Thus, for $|y_{1}|\leq r$, 
\begin{equation*}
|n_{1}-y_{1}|\geq|n_{1}|-r>|n_{1}|/2.
\end{equation*}
Combining this inequality with Proposition \ref{lem:tail}, we obtain
\begin{equation*}
\int_{-r}^{r}k^{[p]}_{t}(n_{1}-y_{1})dy_{1}\leq C^{[p]}_{\delta}t^{\delta/p}\int_{-r}^{r}|n_{1}-y_{1}|^{-1-\delta}dy_{1}\leq 2^{2+\delta}r\ C^{[p]}_{\delta}t^{\delta/p}|n_{1}|^{-1-\delta}.
\end{equation*}
Consequently,
\begin{equation*}
\sum_{|n_{1}|>R}\int_{-r}^{r} k^{[p]}_{t}(n_{1}-y_{1})dy_{1}\leq2^{2+\delta}r\ C^{[p]}_{\delta} t^{\delta/p}\cdot 2\sum_{R<n_{1}\in\mathbb{N}}n_{1}^{-1-\delta}\leq \frac{2^{3+2\delta}r\ C^{[p]}_{\delta}}{\delta}t^{\delta/p}R^{-\delta}.
\end{equation*}
Here, we have used the following integral comparison for real $R>1$:
\begin{equation*}
\sum_{R<j\in\mathbb{N}}j^{-1-\delta}\leq\int_{[R]}^{\infty} x^{-1-\delta}dx
=\frac{[R]^{-\delta}}{\delta}\leq\frac{2^{\delta}}{\delta}R^{-\delta},
\end{equation*}
where $[R]$ denotes the greatest integer not exceeding $R$.\par
On the other hand, Tonelli's theorem is also applicable to the remaining series by nonnegativity. Using the equivalence $\mathbf{1}_{[n_{2}-r,n_{2}+r]}=1$ if and only if $|n_{2}-u|\leq r$, together with $\int k^{[p]}_{t}=1$, we obtain
\begin{align*}
\sum_{n_{2}\in\mathbb{Z}}\int_{-r}^{r} k^{[p]}_{t}(n_{2}-y_{2})\ dy_{2}
&=\sum_{n_{2}\in\mathbb{Z}}\int_{n_{2}-r}^{n_{2}+r}k^{[p]}_{t}(u)du\\
&=\sum_{n_{2}\in\mathbb{Z}}\int_{\mathbb{R}} k^{[p]}_{t}(u)\mathbf{1}_{[n_{2}-r,n_{2}+r]}(u)du\\
&=\int_{\mathbb{R}} k^{[p]}_{t}(u)\sum_{n_{2}\in\mathbb{Z}}\mathbf{1}_{\{|n_{2}-u|\leq r\}}du\\
&=\int_{\mathbb{R}} k^{[p]}_{t}(u)\ \#\{n_{2}\in\mathbb{Z}\mid|n_2-u|\leq r\}du\\
&\leq[2r]+1.
\end{align*}
The case $|n_{2}|>R$ can be treated in exactly the same way. Adding the two estimates yields (\ref{eq:lattice}), with $C_{r,\delta}^{[p]}:=2^{4+2\delta}r C_{\delta}^{[p]}([2r]+1)/\delta$.
\end{proof}
\begin{remark}[Uniform convergence and interchange of limits]
It follows from (\ref{lim:sup}) that, if $\lim_{t\to0+}G_{t}(n)$ exists for each fixed
$n\in\mathbb{Z}^{2}$, then the limit and the lattice point summation can be interchanged: $\lim_{t\to0+}\sum_{n}G_{t}(n)=\sum_{n}\lim_{t\to0+}G_{t}(n)$. Indeed, this follows by decomposing the series into a finite partial sum and a uniformly small tail. Moreover, for any fixed sufficiently large $R>2r+2$, the tail itself tends to zero as $t\to0+$.
\end{remark}

Next, we consider the pointwise limits of $G_{t}(n)$ at lattice points. \par
Let $h=t^{1/p}$. Then, by making the change of variables $z=hy$ in (\ref{eq:scale}), we obtain
\begin{equation}\label{eq:boundary}
G_{t}(n)=\int_{\mathbb{R}^{2}}F_{r}(n-z)K^{[p]}_{t}(z)dz=\int_{\mathbb{R}^{2}}\mathbf{1}_{D_{p}(r)}(n-hy)K^{[p]}_{1}(y)dy,
\end{equation}
namely, $G_{t}$ is obtained by integrating $K^{[p]}_{1}(y)$ only over those $y$ for which
$n-hy$ lies inside the domain. Since $K^{[p]}_{1}\geq0$ and $\int K^{[p]}_{1}=1$, we have $0\leq G_{t}\leq1$. Thus, $G_{t}(n)$ can be interpreted as the proportion of the kernel mass lying inside $D_{p}(r)$, with respect to the weight $K^{[p]}_{1}$, near the lattice point $n$.

\begin{lemma}[\textit{Regularization weights at lattice points}]\label{lem:boundary}\itshape{ }
\begin{equation*}
\lim_{t\to0+}G_{t}(n)=
\begin{cases}
1,&n\in D_{p}(r),\\
0,&n\notin\overline{D_{p}(r)},\\
w_{p}(n),&n\in\partial D_{p}(r),
\end{cases}
\qquad\text{for }n\in\mathbb{Z}^{2},
\end{equation*}
\begin{equation*}
\text{where }w_{p}(n):=
\begin{cases}
1/2, & n_{1}n_{2}\neq 0,\\[2mm]
0, & n_{1}n_{2}=0,\quad 0<p<1,
\end{cases}
\quad
\begin{cases}
1/4, & n_{1}n_{2}=0,\quad p=1,\\[2mm]
1/2, & n_1n_2=0,\quad p=2.
\end{cases}
\end{equation*}
\end{lemma}
\begin{proof}
(i) \textit{The interior case}. Suppose that $n\in D_{p}(r)$. By the openness of $D_{p}(r)$, there exists $\varepsilon>0$ such that $B(n,\varepsilon)\subset D_{p}(r)$. If $|z|_{2}<\varepsilon$, then $|(n-z)-n|_{2}=|z|_{2}<\varepsilon$, and hence $n-z\in D_{p}(r)$. Therefore, by (\ref{eq:approx}), we have
\begin{equation*}
1\geq G_{t}(n)\geq\int_{|z|_2<\varepsilon}K^{[p]}_{t}(z)dz=1-\int_{|z|_2\geq\varepsilon}K^{[p]}_{t}(z)dz\longrightarrow1\qquad\text{as }t\to0+.
\end{equation*}
(ii) \textit{The exterior case}. Suppose that $n\notin\overline{D_{p}(r)}$. Then there exists $\varepsilon>0$ such that $B(n,\varepsilon)\cap D_{p}(r)=\emptyset$. By an argument similar to that in (i), we have $n-z\notin D_{p}(r)$ whenever $|z|_{2}<\varepsilon$. Hence, by (\ref{eq:approx}), we obtain
\begin{equation*}
0\leq G_{t}(n)=\int_{|z|_{2}\geq\varepsilon}F_{r}(n-z)K^{[p]}_{t}(z)dz\leq\int_{|z|_{2}\geq\varepsilon}K^{[p]}_{t}(z)dz\longrightarrow0\qquad\text{as }t\to0+.
\end{equation*}
(iii) \textit{The boundary case}. First, suppose that $n$ is not on either coordinate axis, namely, $n_{1}n_{2}\neq0$. Let $\Psi(x)=|x_{1}|^{p}+|x_{2}|^{p}-r^{p}$ (so that $D_{p}(r)=\{x\in\mathbb{R}^{2}\mid\Psi(x)<0\}$). Then $\Psi(n)=0$ and $\nabla\Psi(n)=(p\ \mathrm{sgn}(n_{1})|n_{1}|^{p-1},p\ \mathrm{sgn}(n_{2})|n_{2}|^{p-1})\neq0$. For each fixed $y$, we have
\begin{equation*}
\frac{\Psi(n-hy)}h=-\nabla\Psi(n)\cdot y+o(1)\qquad\text{as }h\to0+.
\end{equation*}
If $c(y):=\nabla\Psi(n)\cdot y>0$, then the right-hand side converges to $-c(y)<0$, so $\Psi(n-hy)$ is negative for sufficiently small $h>0$. Similarly, if $c(y)<0$, then $\Psi(n-hy)$ is positive for sufficiently small $h>0$. The remaining exceptional set $c(y)=0$ is a straight line and hence has two-dimensional Lebesgue measure zero. In this case, the integrand in (\ref{eq:boundary}) is dominated by $K^{[p]}_{1}\in L^{1}(\mathbb{R}^{2})$, and the following limit holds. Thus, the dominated convergence theorem can be applied.
\begin{equation*}
\mathbf{1}_{D_{p}(r)}(n-hy)\longrightarrow\mathbf{1}_{\{\nabla\Psi(n)\cdot y>0\}}(y)\qquad\text{as }h\to0+.
\end{equation*}
Moreover, since $K^{[p]}_{1}(-y)=K^{[p]}_{1}(y)$, the integrals over the two half-planes separated by the line $\nabla\Psi(n)\cdot y=0$, which passes through the origin, are equal. Hence, we obtain
\begin{equation*}
\lim_{t\to0+}G_{t}(n)=\int_{\{\nabla\Psi(n)\cdot y>0\}}K^{[p]}_{1}(y)dy=\frac{1}{2}\int_{\mathbb{R}^{2}}K^{[p]}_{1}(y)dy=\frac{1}{2}.
\end{equation*}
\par
We next consider points on the coordinate axes ($n_{1}n_{2}=0$). By symmetry, it suffices to consider $n=(r,0)$. When $p=2$, we have $\nabla\Psi(r,0)=(2r,0)\neq0$, and hence the same half-plane argument as for a smooth boundary point yields $\lim_{t\to0+}G_{t}(t)=1/2$. \par
On the other hand, for $0<p\leq1$, the condition that $n-hy$ lies inside $D_{p}(r)$ is given by
\begin{equation}\label{condition}
|r-hy_{1}|^{p}+h^{p}|y_{2}|^{p}<r^{p}.
\end{equation}
For $0<p<1$ (where the boundary point on the coordinate axis is a cusp), by the Taylor expansion for fixed $y$,
\begin{equation*}
|r-hy_{1}|^{p}=r^{p}|1-(y_{1}/r)h|^{p}=r^{p}-pr^{p-1}hy_{1}+O_{y}(h^{2})\qquad\text{as }h\to0+
\end{equation*}
and dividing by $h>0$, the interior condition (\ref{condition}) reduces to the following.
\begin{equation*}
-pr^{p-1}y_{1}+h^{p-1}|y_{2}|^{p}+O_{y}(h)<0.
\end{equation*}
In this case, if $y_{2}\neq0$, then $h^{p-1}|y_{2}|^{p}\to+\infty$, and hence the indicator function of the interior of $D_{p}(r)$ converges to zero. That is,
\begin{equation*}
\mathbf{1}_{D_{p}(r)}(n-hy)\longrightarrow0\qquad\text{for }y_{2}\neq0,\ \text{as }h\to0+.
\end{equation*}
Since the exceptional set $\{y_{2}=0\}$ is a straight line and has two-dimensional Lebesgue measure zero, the dominated convergence theorem yields $\lim_{t\to0+}G_{t}(r,0)=0$.\par\vspace{5pt}
On the other hand, consider the case $p=1$ (where the boundary point on the coordinate axis is a vertex of the diamond). At $n=(r,0)$, we have $r-hy_{1}>0$ for sufficiently small $h>0$, for each fixed $y$. Thus, the interior condition (\ref{condition}) becomes
\begin{equation*}
r-hy_{1}+h|y_{2}|<r\quad\Longleftrightarrow\quad y_{1}>|y_{2}|.
\end{equation*}
The four sets
\begin{equation*}
H_{1}:=\{y_{1}>|y_{2}|\},\quad H_{2}:=\{-y_{1}>|y_{2}|\},\quad H_{3}:=\{y_{2}>|y_{1}|\},\quad H_{4}:=\{-y_{2}>|y_{1}|\}
\end{equation*}
partition the plane into four regions, apart from their boundary lines. Since $K^{[1]}_{1}(y_{1},y_{2})=k^{[1]}_{1}(y_{1})k^{[1]}_{1}(y_{2})$ is invariant under coordinate interchange and sign changes, the integrals over these four regions are equal. Consequently, we obtain the desired conclusion. 
\begin{equation*}
\lim_{t\to0+}G_{t}(n)=\int_{H_{1}}K^{[p]}_{1}(y)dy=\frac{1}{4}\int_{\mathbb{R}^{2}}K^{[p]}_{1}(y)dy=\frac{1}{4}.
\end{equation*}
\end{proof}

\section{Proof of Theorem \ref{thm:abel-GHI}}
\hspace{13pt}We begin the proof of our main theorem (Theorem \ref{thm:abel-GHI}) by establishing the regularized Poisson summation identity, which constitutes a fundamental ingredient in its proof. The underlying principle is the classical Poisson summation formula applied after convolution with an anisotropic regularization kernel. Since the kernel and the domain under consideration have specific analytic and geometric features, we give a direct proof in order to verify the absolute convergence of the relevant series and justify the interchange of summation and integration. The resulting identity provides the starting point for the subsequent analysis of the regularization limit and the boundary contributions.
\begin{lemma}[\textit{Regularized Poisson Summation Identity}]\label{lem:poisson}
\itshape{For every $t>0$, the following identity holds, with both sides converging absolutely:}
\begin{equation}\label{eq:poisson}
\sum_{n\in\mathbb{Z}^{2}}(F_{r}\ast K^{[p]}_{t})(n)=\sum_{m\in\mathbb{Z}^{2}}\widehat{F_{r}}(m)e^{-t|m|^{p}_{p}}.
\end{equation}
\end{lemma}
\begin{proof}
Put $G_{t}:=F_{r}\ast K^{[p]}_{t}$. For each $y$, if $x_{i}\to x$, then $K^{[p]}_{t}(x_{i}-y)\to K^{[p]}_{t}(x-y)$ holds. Moreover, we have $|F_{r}(y)K^{[p]}_{t}(x_{i}-y)|\leq \sup_{x\in\mathbb{R}^{2}}|K^{[p]}_{t}(x)|\cdot F_{r}(y)$ and the right-hand side is integrable. Thus, the dominated convergence theorem implies that $G_{t}$ is continuous. \par
Next, fix a compact set $C\subset\mathbb{R}^{2}$, and let $M>0$ be such that $|x_{1}|,|x_{2}|\leq M$ for all $x\in C$. If $y\in D_{p}(r)$, then $|y_{1}|,|y_{2}|\leq r$. In this case, if $|n_{1}|,|n_{2}|\geq2(M+r)+1$, then
\begin{equation*}
|x_{j}+n_{j}-y_{j}|\geq|n_{j}|-|x_{j}|-|y_{j}|\geq|n_{j}|-(M+r)\geq\frac{1}{2}(|n_{j}|+1)\quad\text{for }j=1,2.
\end{equation*}
Hence, by Proposition \ref{lem:tail}, we obtain
\begin{equation}\label{inq-n}
k^{[p]}_{t}(x_{j}+n_{j}-y_{j})\leq2^{1+\delta}C_{\delta}^{[p]}t^{\delta/p}(1+|n_{j}|)^{-1-\delta}\quad\text{for }j=1,2.
\end{equation}
On the other hand, recalling that $k^{[p]}_{t}(x) \leq(2/p)\Gamma(1/p)t^{-1/p}$, we see that, by enlarging the constant if necessary, the estimate (\ref{inq-n}) also holds for the finitely many remaining values of $n_{j}$ for which $|n_{j}|<2(M+r)+1$, for each $j=1,2$. Consequently,
\begin{align}
\sup_{x\in C}|G_{t}(x+n)|&=\sup_{x\in C}\int_{D_{p}(r)}k^{[p]}_{t}(x_{1}+n_{1}-y_{1})k^{[p]}_{t}(x_{2}+n_{2}-y_{2})dy\notag\\
&\leq C_{M,r,t,\delta}^{[p]}(1+|n_1|)^{-1-\delta}(1+|n_{2}|)^{-1-\delta}.\label{eq:periodbound}
\end{align}
By the integral test, we have $\sum_{j\in\mathbb{Z}}(1+|j|)^{-1-\delta}\leq1+2/\delta$. It therefore follows from (\ref{eq:periodbound}) that $\sum_{n\in\mathbb{Z}^{2}}\sup_{x\in C}|G_{t}(x+n)|\leq(1+2/\delta)^{2}C_{C,r,t,\delta}^{[p]}<+\infty$. Thus, by the Weierstrass $M$-test, the periodization
\begin{equation*}
g_{t}(x)=\sum_{n\in\mathbb{Z}^{2}}G_{t}(x+n)
\end{equation*}
converges uniformly on compact subsets of $\mathbb{R}^{2}$. In particular, $g_{t}$ is continuous and $1$-periodic. For $m\in\mathbb{Z}^{2}$, the absolute convergence of the periodization allows us to interchange the summation and integration when computing its Fourier coefficients. Then, using $e^{-2\pi im\cdot n}=1$ and the change of variables $u=x+n$, we obtain
\begin{align*}
\widehat{g_{t}}(m)&=\int_{[0,1]^2}\sum_{n\in\mathbb{Z}^{2}}G_{t}(x+n)e^{-2\pi im\cdot x}dx\\
&=\sum_{n\in\mathbb{Z}^{2}}\int_{[0,1]^2+n}G_{t}(u)e^{-2\pi im\cdot u}du=\widehat{G_{t}}(m)\\
&=\widehat{(F_{r}\ast K^{[p]}_{t})}(m)=\widehat{F_{r}}(m)\widehat{K^{[p]}_{t}}(m)=\widehat{F_{r}}(m)e^{-t|m|^{p}_{p}}.
\end{align*}
The convolution theorem for the Fourier transform is justified by Fubini's theorem, since $F_{r}, K^{[p]}_{t}\in L^{1}(\mathbb{R}^{2})$. For the Fourier transform of $K^{[p]}_{t}$, see Section 2. Moreover, since $|\widehat{F_{r}}(m)|\leq\int|F_{r}|=A_{p}r^{2}$, where $A_{p}:=2\Gamma(1/p)^{2}/(p\Gamma(2/p))$ is the area of the unit $p$-circle, we have
\begin{align*}
\sum_{m\in\mathbb{Z}^{2}}|\widehat{g_{t}}(m)|&\leq A_{p}r^{2}\biggl(\sum_{j\in\mathbb{Z}}e^{-t|j|^{p}}\biggr)^{2}\\
&\leq A_{p}r^{2}\biggl(1+2\int_{0}^{\infty}e^{-tx^{p}}dx\biggr)^{2}
=A_{p}r^{2}\biggl(1+\frac{2}{p}\Gamma(1/p)t^{-1/p}\biggr)^{2}<+\infty.
\end{align*}
Consequently, the series $h_{t}(x):=\sum_{m}\widehat{g_{t}}(m)e^{2\pi im\cdot x}$ converges absolutely and uniformly, and hence defines a continuous periodic function. By the orthogonality of the exponential functions, we obtain
\begin{equation*}
\widehat{h_{t}}(\ell)=\int_{[0,1]^{2}}\biggl(\sum_{m\in\mathbb{Z}^{2}}\widehat{g_{t}}(m)e^{2\pi im\cdot x}\biggr)e^{-2\pi i\ell\cdot x}dx=\sum_{m\in\mathbb{Z}^{2}}\widehat{g_{t}}(m)\int_{[0,1]^{2}}e^{2\pi i(m-\ell)\cdot x}dx=\widehat{g_{t}}(\ell).
\end{equation*}
By the uniqueness of Fourier coefficients for continuous periodic functions, we conclude that $h_{t}=g_{t}$. Setting $x=0$, we obtain (\ref{eq:poisson}).
\end{proof}
\begin{remark}
The regularized Poisson summation identity above is not specific to the indicator function of a $p$-circle domain. The same argument applies to any bounded, compactly supported measurable function $f$ on $\mathbb{R}^{2}$, with $F_r$ replaced by $f$. Thus, the identity itself is a particular instance of the classical periodization principle in Fourier analysis.
The geometric features specific to the present problem arise in the limit as the regularization parameter tends to zero. In particular, the boundary behavior of the $p$-circle determines the limiting contributions of lattice points, which differ between smooth boundary points and axial cusps. The rigorous identification of these contributions, together with the $p$-Bessel Fourier representation, is essential for establishing the Abel-summed GHI with explicit boundary corrections.
\end{remark}

\begin{proof}[Proof of Theorem \ref{thm:abel-GHI}]
We first consider the limit as $t\to0+$ of the left-hand side of Lemma \ref{lem:poisson} (Regularized Poisson Summation Identity). Fix $R>2r+2$. Then all lattice points belonging to $D_{p}(r)$ or its boundary are contained in the region $|n|_{\infty}\leq R$. By Lemma \ref{lem:boundary}, we have, for the finite sum,
\begin{equation*}
\lim_{t\to0+}\sum_{|n|_{\infty}\leq R}G_{t}(n)=\sum_{|n|_{\infty}\leq R}\lim_{t\to0+}G_{t}(n)=N_{p}(r)+B_{p}(r).
\end{equation*}
On the other hand, Lemma \ref{lem:lattice} yields
\begin{equation*}
0\leq\sum_{|n|_{\infty}>R}G_{t}(n)\leq C_{r,\delta}^{[p]}t^{\delta/p}R^{-\delta}\longrightarrow0\qquad\text{as }t\to0+.
\end{equation*}
Consequently,
\begin{equation*}\label{eq:limsum}
\lim_{t\to0+}\sum_{n\in\mathbb{Z}^{2}}G_{t}(n)=N_{p}(r)+B_{p}(r).
\end{equation*}
Since the $m=0$ term on the right-hand side of the Regularized Poisson Summation Identity is given by $\widehat{F_{r}}(0)=A_{p}r^{2}$, this identity can be rewritten, by the definition of $P_{p}$, as
\begin{equation}\label{eq:sumFT}
P_{p}(r)=\lim_{t\to0+}\sum_{m\in\mathbb{Z}^{2}\setminus\{0\}}\widehat{F_{r}}(m)e^{-t|m|^{p}_{p}}-B_{p}(r).
\end{equation}
We now employ the $p$-Bessel-Fourier transform representation of $p$-radial functions given in \cite[(2.3)]{K3}. For $m\neq0$, putting $|m|^{p}_{p}=s$ and denoting the corresponding distorted angle by $\varphi_{m}$, the Fourier transform of $F_{r}$ can be expressed as follows.
\begin{align*}
    \widehat{F_{r}}(m)&=p\Gamma(1/p)^{2}\int_{0}^{\infty}\mathcal{J}_{0,\varphi_{m}}^{[p]}(2\pi \rho s^{\frac{1}{p}})\phi_{F_{r}}(\rho)\rho\ d\rho\\
    &=p\Gamma(1/p)^{2}\int_{0}^{r}\mathcal{J}_{0,\varphi_{m}}^{[p]}(2\pi \rho s^{\frac{1}{p}})\rho\ d\rho\\
    &=p\Gamma(1/p)^{2}\int_{0}^{r}(2\pi \rho s^{\frac{1}{p}})\mathcal{J}_{0,\varphi_{m}}^{[p]}(2\pi \rho s^{\frac{1}{p}})d\rho\cdot (2\pi s^{1/p})^{-1}\\
    &=p\Gamma(1/p)^{2}\int_{0}^{2\pi rs^{1/p}}\tau\mathcal{J}_{0,\varphi_{m}}^{[p]}(\tau)d\tau\cdot (2\pi s^{1/p})^{-2}\\
    &=\frac{p\Gamma(1/p)^{2}}{(2\pi s^{1/p})^{2}}(2\pi rs^{\frac{1}{p}})\mathcal{J}_{1,\varphi_{m}}^{[p]}(2\pi rs^{\frac{1}{p}})\\
    &=\frac{p\Gamma(1/p)^{2}}{2\pi}\frac{r}{s^{1/p}}\mathcal{J}_{1,\varphi_{m}}^{[p]}(2\pi rs^{\frac{1}{p}}).
  \end{align*}
In the passage from the fourth to the fifth line, we used the integral form of the differentiation formula for $p$-Bessel functions \cite[Proposition 2.6]{K3}, namely, $\int_{0}^{r}\tau\mathcal{J}_{0,\varphi}^{[p]}(\tau)d\tau=r\mathcal{J}_{1,\varphi}^{[p]}(r)$.
For each $t>0$, we have
\begin{equation*}
\sum_{m\neq0}|\widehat{F_{r}}(m)|e^{-t|m|^{p}_{p}}\leq A_{p}r^{2}\sum_{m\neq0}e^{-t|m|^{p}_{p}}<+\infty.
\end{equation*}
Thus, the absolute convergence justifies regrouping the summation over $m$ according to $s=|m_{1}|^{p}+|m_{2}|^{p}$. Moreover, the lattice points corresponding to the same $s$ can be indexed by the distorted angles
\begin{equation*}\vspace{-5pt}
  \varphi\in\mathcal{A}_{s}^{[p]}=\left\{\varphi\in[0,2\pi)\ \middle|\ \left(\mathrm{sgn}\left(\cos\varphi\right)s^{\frac{1}{p}}|\cos\varphi|^{\frac{2}{p}},\mathrm{sgn}\left(\sin\varphi\right)s^{\frac{1}{p}}|\sin\varphi|^{\frac{2}{p}}\right)\in\mathbb{Z}^{2}\right\}.
\end{equation*}
Let $\mu$ denote the counting measure on the set of $s$ satisfying $\mathcal{A}^{[p]}_{s}\neq\emptyset$, assigning mass $1$ to each such $s$. Then
\begin{equation*}
\sum_{m\neq0}\widehat{F_{r}}(m)e^{-t|m|_{p}^{p}}=\frac{p\Gamma(1/p)^2}{2\pi}r\int_{1}^{\infty}\frac{e^{-ts}}{s^{1/p}}\biggl(\sum\nolimits_{\varphi\in\mathcal{A}_{s}^{[p]}}\mathcal{J}_{1,\varphi}^{[p]}(2\pi rs^{\frac{1}{p}})\biggr)d\mu(s).
\end{equation*}
Substituting this into (\ref{eq:sumFT}), we obtain the regularized GHI (\ref{eq:abel-GHI}), which is the main result of this paper.
\end{proof}
\begin{remark}[Explicit Form of the Boundary Correction]
For $0<p<1$, we have $B_{p}(r)=\frac{1}{2}\#\{n\in\partial D_{p}(r)\cap\mathbb{Z}^{2}\mid n_{1}n_{2}\neq0\}$. For $p=1$, $B_{p}(r)$ is given by the sum of one-half of the number of non-axis boundary lattice points and one-quarter of the number of axial vertices. For $p=2$, $B_{p}(r)$ is equal to one-half of the total number of boundary lattice points. If there are no lattice points on the boundary, then $B_{p}(r)=0$.
\end{remark}

\section{Future Directions}
\hspace{13pt}A natural question left open by the present work is whether the Abel-summed generalized Hardy-type identity can be strengthened to an identity with ordinary convergence of the unregularized oscillatory series, at least when the boundary of the corresponding $p$-circle contains no lattice points. The regularized identity established here does not by itself imply convergence of sharp $p$-radial Fourier partial sums. Establishing such convergence would require additional control of oscillatory cancellation beyond the convergence of the regularized series.\par
The anisotropic geometry of astroid-type p-circles presents a further challenge. For
$0<p<1$, their boundaries have cusps on the coordinate axes. This anisotropy is reflected in the integral representation of the Fourier transform of $p$-radial functions, whose kernel involves the distorted-angle-dependent $p$-Bessel function of order zero. 
  \begin{equation*}
    \hat{F}(x)=p\Gamma(1/p)^{2}\int_{0}^{\infty}\mathcal{J}_{0,\varphi}^{[p]}(2\pi \rho r)\phi_{F}(\rho)\rho\ d\rho 
  \end{equation*}
for $x=(\mathrm{sgn}(\cos\varphi)r|\cos\varphi|^{\frac{2}{p}},\mathrm{sgn}(\sin\varphi)r|\sin\varphi|^{\frac{2}{p}})$. More concretely, the associated oscillatory integral representations have stationary points whose locations depend on the direction of the Fourier variable; as this direction approaches a coordinate axis, an interior stationary point may approach an endpoint of integration \cite[Section 2.2, Table 1]{K2}. For astroid-type $p$-circles, uniform estimates for the $p$-Bessel function of order zero, including axial directions, have already been established by analyzing higher-order derivatives of the phase and applying \textit{van der Corput's lemma}\cite[Proposition 2.8 and Theorem 1.5]{K2}. Together with the previously obtained axial asymptotic formula\cite[Theorem 1.3]{K4}
\begin{align*}
  \mathcal{J}_{\omega,\varphi_{\mathrm{axis}}}^{[p]}(r)=\frac{4\Gamma(p+1)}{p^{\omega+1}\Gamma(\omega+\frac{1}{p}-1)\Gamma(\frac{1}{p})}&r^{\omega-p-1}\sin\frac{\pi p}{2}+\frac{4p^{\frac{1}{p}-2}}{\Gamma(\frac{1}{p})}r^{-\frac{1}{p}}\cos\Bigl(r-\frac{\pi}{2}(\frac{1}{p}+\omega)\Bigr)\\
  &+\mathcal{O}(r^{-\beta(p,\omega)})\quad\text{for }\omega\geq0,\ \varphi_{\mathrm{axis}}\in\left\{0,\frac{\pi}{2},\pi,\frac{3\pi}{2}\right\},
\end{align*}
where $\beta(p,\omega):=\min\{1+2p-\omega,1+1/p\}$, these results motivate further uniform asymptotic analysis of the first-order function across axial and non-axial directions. Such analysis, combined with suitable estimates for oscillatory lattice-point sums, may help address ordinary convergence of the generalized Hardy-type identity and related localization or \textit{Pinsky}-type questions for two-dimensional $p$-radial Fourier partial sums.

\vspace{25pt}\hspace{-15pt}\textbf{Data availability.}\quad No data was used for the research described in the article.

\itshape{
  \hspace{13pt}The author's affiliation: Graduate School of Mathematics, Nagoya University, Chikusa-ku, Nagoya 464-8602, Japan\\
  \hspace{13pt}The author's email address: kitajima.masaya.p2@a.mail.nagoya-u.ac.jp}
\end{document}